%% file: SC.tex
\documentclass[reqno]{amsart}

\pdfoutput=1
\usepackage{general,resistance,ifthen,pxfonts,leading}
\usepackage[protrusion=true,expansion]{microtype}
\usepackage{cancel,xy}
\usepackage[switch*,pagewise]{lineno}

\usepackage[bookmarks, colorlinks=true, pdfstartview=FitV, linkcolor=blue, citecolor=blue, urlcolor=blue]{hyperref}
\usepackage{cite}  

\usepackage[left=1.85in,right=1.85in,top=1.67in,bottom=1.40in]{geometry}

\numberwithin{equation}{section} \numberwithin{theorem}{section}

\newcommand{\Lb}{{\ensuremath{\Lap^{(b)}}}\xspace}     
\newcommand{\Lc}{{\ensuremath{\Lap^{(c)}}}\xspace}     
\newcommand{\Lbc}{{\ensuremath{\Lap^{(b)^{\bigstar_{\negsp[2]c}}}}}\xspace}     
\newcommand{\Eb}{{\ensuremath{\energy^{(b)}}}\xspace}  
\newcommand{\Ec}{{\ensuremath{\energy^{(c)}}}\xspace}  
\newcommand{\Finb}{{\ensuremath{\Fin^{(b)}}}\xspace}  
\newcommand{\Finc}{{\ensuremath{\Fin^{(c)}}}\xspace}  
\newcommand{\Harmb}{{\ensuremath{\Harm^{(b)}}}\xspace}  
\newcommand{\Harmc}{{\ensuremath{\Harm^{(c)}}}\xspace}  
\newcommand{\HEb}{\ensuremath{\sH_\Eb}\xspace}         
\newcommand{\HEc}{\ensuremath{\sH_\Ec}\xspace}         
\newcommand{\vxb}[1][x]{\ensuremath{v_{#1}^{(b)}}\xspace}       
\newcommand{\vxc}[1][x]{\ensuremath{v_{#1}^{(c)}}\xspace}       
\newcommand{\vyb}[1][y]{\ensuremath{v_{#1}^{(b)}}\xspace}       
\newcommand{\vyc}[1][y]{\ensuremath{v_{#1}^{(c)}}\xspace}       
\newcommand{\gmb}{\ensuremath{\gm^{(b)}}\xspace}       
\newcommand{\gmc}{\ensuremath{\gm^{(c)}}\xspace}       
\newcommand{\Hom}{\ensuremath{\operatorname{Hom}}\xspace}            
\newcommand{\sgn}{\ensuremath{\operatorname{sgn}}\xspace} 

\begin{document}
  \leading{14pt} 
  \linenumbers


\title[Spectral comparisons of networks]{Spectral comparisons between networks with different conductance functions}

\author{Palle E. T. Jorgensen {\protect \and} Erin P. J. Pearse}

\address{Palle E. T. Jorgensen, University of Iowa, Iowa City, IA 52246-1419 USA}
\email{palle-jorgensen@uiowa.edu}

\address{Erin P. J. Pearse, California Polytechnic University, San Luis Obispo, CA 93407-0403 USA}
\email{epearse@calpoly.edu}

\thanks{The work of PETJ was partially supported by NSF grant DMS-0457581. The work of EPJP was partially supported by the University of Iowa Department of Mathematics NSF VIGRE grant DMS-0602242.}

\begin{abstract}
  For an infinite network consisting of a graph with edge weights prescribed by a given conductance function $c$, we consider the effects of replacing these weights with a new function $b$ that satisfies $b \leq c$ on each edge. In particular, we compare the corresponding energy spaces and the spectra of the Laplace operators acting on these spaces. We use these results to derive estimates for effective resistance on the two networks, and to compute a spectral invariant for the canonical embedding of one energy space into the other.
\end{abstract}



\date{\bf\today.}

\maketitle

\setcounter{tocdepth}{1}
{\small \tableofcontents}

\allowdisplaybreaks

\input{introduction}

\input{setup}

\input{conductances}

\input{permanence}

\input{examples}

\subsection*{Acknowledgements}

The authors are grateful for stimulating comments, helpful advice, and valuable references from Dorin Dutkay, Paul Muhly, and others. We also thank the students and colleagues who have endured our talks on this material and raised fruitful questions.

\bibliographystyle{alpha}
\bibliography{networks}

\label{lastpage}
\end{document}

%% file: introduction.tex

\section{Introduction}
\label{sec:intro}

We begin with a network structure defined by a set of vertices $G$ and a conductance function $c:G \times G \to \bR^+$ which specifies the both the adjacency relation and the edge weights; two vertices $x$ and $y$ are neighbours iff $c_{xy}>0$. The case of  primary interest is when $G$ is infinite, in which case the Hilbert space \HE (comprised of functions of finite Dirichlet energy) has a rich structure and the Laplace operator \Lap corresponding to the network may be an unbounded operator on \HE. (Precise definitions of these terms may be found in  Definition~\ref{def:graph-laplacian}, Definition~\ref{def:graph-energy}, and Definition~\ref{def:H_energy}.)


The Hilbert space \HE has a rather different geometry than the more familiar $\ell^2(G)$, and depends crucially on the choice of conductance function $c$. The same is true for the Laplacian \Lap as a linear operator on \HE. 
In this paper, we use the framework developed in earlier projects (see \cite{DGG,ERM,bdG,SRAMO,Friedrichs,Multipliers,LPS, Interpolation, RBIN, OTERN}) to compute certain spectral theoretic information; as well as resistance metrics on the underlying vertex set. In particular, we explore how certain quantities depend on the choice of $c$, in comparison to another conductance function, which we denote by $b$. It will be assumed that both $b$ and $c$ yield a connected weighted graph, although we allow for the case when $c_{xy}>0$ and $b_{xy}=0$ (so that $x$ and $y$ are neighbours in $(G,c)$ but not in $(G,b)$).
The data, defined from $b$ and $c$, to be compared are as follows:
\begin{enumerate}
  \item the energy forms \Eb and \Ec, and the respective energy Hilbert spaces \HEb and \HEc that they define;
  \item the systems of dipole vectors that form reproducing kernels for the two Hilbert spaces; see Definition~\ref{def:energy-kernel};
  \item the respective Laplace operators \Lb and \Lc, and their spectra;
  \item the spaces of finite-energy harmonic functions on \HEb and \HEc; and  
  \item the effective resistance metrics on \HEb and \HEc.
\end{enumerate}
We focus our study on the case when one of the two energy-Hilbert spaces is contractively contained in the other, which corresponds to the inequality $b \leq c$. In this case, we believe that our results have applications to percolation theory and the study of random walks in random environments, as well as to dilation theory (see \cite{Arveson-dilation}) and the contractive inclusion of Hilbert spaces (see \cite{Sarason94}).
     
Of special operator theoretic significance is the adjoint of the contractive inclusion mapping. The issues involved with the adjoint operator are subtle as the computation of the adjoint operator depends on the Hilbert-inner products used. It is the adjoint operator that allows one to compute the respective systems of dipole vectors that form reproducing kernels for the two Hilbert spaces; see Definition~\ref{def:energy-kernel}.
We further derive an invariant (involving induced linear maps between the respective spaces of finite-energy harmonic functions) which distinguishes two networks when $G$ is fixed and the conductance functions vary.

We also give a necessary and sufficient condition on a fixed conductance function $c$ having its energy Hilbert space \Ec boundedly contained in \HEb ($b = 1$); i.e., contractive containment in the ``flat'' energy Hilbert space corresponding to constant conductance $b$. The significance of this is that computations in \HEb are typically much easier, and that the conclusions obtained there may then be transferred to \HEc.    

Our results are illustrated with concrete examples.

%% file: setup.tex

\section{Basic terms and previous results}
\label{sec:Basic-terms-and-previous-results}

In this section, we introduce the key notions used throughout this paper: resistance networks, the energy form \energy, the Laplace operator \Lap, and their elementary properties. 

\begin{defn}\label{def:ERN}
  A \emph{(resistance) network} is a connected graph $(\Graph,\cond)$, where \Graph is a graph with vertex set \verts, and \cond is the \emph{conductance function} which defines adjacency by $x \nbr y$ iff $c_{xy}>0$, for $x,y \in \verts $. We assume $\cond_{xy} = \cond_{yx} \in [0,\iy)$, and write $\cond(x) := \sum_{y \nbr x} \cond_{xy}$. We require that the graph is \emph{locally finite}, i.e., that every vertex has only finitely many neighbors. 
\end{defn}

In this definition, connected means simply that for any $x,y \in \verts $, there is a finite sequence $\{x_i\}_{i=0}^n$ with $x=x_0$, $y=x_n$, and $\cond_{x_{i-1} x_i} > 0$, $i=1,\dots,n$.  
We may assume there is at most one edge from $x$ to $y$, as two conductors $\cond^1_{xy}$ and $\cond^2_{xy}$ connected in parallel can be replaced by a single conductor with conductance $\cond_{xy} = \cond^1_{xy} + \cond^2_{xy}$. Also, we assume $\cond_{xx}=0$ so that no vertex has a loop. 

Since the edge data of $(\Graph,\cond)$ is carried by the conductance function, we will henceforth simplify notation and write $x \in \Graph$ to indicate that $x$ is a vertex. For any network, one can fix a reference vertex, which we shall denote by $o$ (for ``origin''). It will always be apparent that our calculations depend in no way on the choice of $o$.

\begin{defn}\label{def:graph-laplacian}
  The \emph{Laplacian} on \Graph is the linear difference operator which acts on a function $v:\Graph \to \bC$ by
  \linenopax
  \begin{equation}\label{eqn:Lap}
    (\Lap v)(x) :
    = \sum_{y \nbr x} \cond_{xy}(v(x)-v(y)).
  \end{equation}
  The domain of \Lap is discussed in detail in \eqref{eqn:Lap-domains}, below. A function $v:\Graph \to \bC$ is \emph{harmonic} iff $\Lap v(x)=0$ for each $x \in \Graph$. 
\end{defn}

It is clear from \eqref{eqn:Lap} that \Lap commutes with complex conjugation and therefore the subspace of real-valued functions is invariant under \Lap.

We have adopted the physics convention (so that the spectrum is nonnegative) and thus our Laplacian is the negative of the one commonly found in the PDE literature. The network Laplacian \eqref{eqn:Lap} should not be confused with the stochastically renormalized Laplace operator $\cond^{-1} \Lap$ which appears in the probability literature, or with the spectrally renormalized Laplace operator $\cond^{-1/2} \Lap \cond^{-1/2}$ which appears in the literature on spectral graph theory (e.g., \cite{Chung}).

\begin{defn}\label{def:graph-energy}
  The \emph{energy} of functions $u,v:\Graph \to \bC$ is given by the (closed, bilinear) Dirichlet form
  \linenopax
  \begin{align}\label{eqn:def:energy-form}
    \energy(u,v)
    := \frac12 \sum_{x,y \in \Graph} \cond_{xy}(\cj{u}(x)-\cj{u}(y))(v(x)-v(y)),
  \end{align}
  with the energy of $u$ given by $\energy(u) := \energy(u,u)$.
  The \emph{domain} of the energy form is
  \linenopax
  \begin{equation}\label{eqn:def:energy-domain}
    \dom \energy = \{u:\Graph \to \bC \suth \energy(u)<\iy\}.
  \end{equation}
\end{defn}
Note that \eqref{eqn:def:energy-form} converges iff it converges absolutely, by the Schwarz inequality, so this summation is well defined. Since $\cond_{xy}=\cond_{yx}$ and $\cond_{xy}=0$ for nonadjacent vertices, the initial factor of $\frac12$ in \eqref{eqn:def:energy-form} implies there is exactly one term in the sum for each edge in the network. 

\subsection{The energy space \HE} 
\label{sec:The-energy-space}

The energy form \energy is sesquilinear and conjugate symmetric on $\dom \energy$ and would be an inner product if it were positive definite.   

\begin{defn}\label{def:H_energy}\label{def:The-energy-Hilbert-space}
  Let \one denote the constant function with value 1 and recall that $\ker \energy = \bC \one$. Then $\HE := \dom \energy / \bC \one$ is a Hilbert space with inner product and corresponding norm given by
  \linenopax
  \begin{equation}\label{eqn:energy-inner-product}
    \la u, v \ra_\energy := \energy(u,v)
    \q\text{and}\q
    \|u\|_\energy := \energy(u,u)^{1/2}.
  \end{equation}
  We call \HE the \emph{energy (Hilbert) space}. 
\end{defn}

\begin{remark}\label{rem:completion}
  Since \Graph is connected, it is possible to show (with the use of Fatou's lemma) that $\dom \energy / \bC \one$ is complete; 
  see \cite{DGG,OTERN} for further details regarding this point.
\end{remark}

\begin{defn}\label{def:vx}\label{def:energy-kernel}
  Let $v_x$ be defined to be the unique element of \HE for which
  \linenopax
  \begin{equation}\label{eqn:v_x}
    \la v_x, u\ra_\energy = u(x)-u(o),
    \qq \text{for every } u \in \HE.
  \end{equation}
  The existence and uniqueness of $v_x$ for each $x \in G$ is implied by the Riesz lemma. It follows from \eqref{eqn:v_x} that $\{v_x\}_{x \in \Graph}$ forms a reproducing kernel for \HE (called the \emph{energy kernel}; see \cite[Cor.~2.7]{DGG}) and that $\spn\{v_x\}_{x \in \Graph}$ is dense in \HE. 
\end{defn}

Note that $v_o$ corresponds to a constant function, since $\la v_o, u\ra_\energy = 0$ for every $u \in \HE$. Therefore, $v_o$ may often be safely ignored or omitted during calculations. 

\begin{defn}\label{def:dipole}
  A \emph{dipole} is any $v \in \HE$ satisfying the pointwise identity $\Lap v = \gd_x - \gd_y$ for some vertices $x,y \in \Graph$. One can check that $\Lap v_x = \gd_x - \gd_o$; cf. \cite[Lemma~2.13]{DGG}.
\end{defn}

Note that dipoles always exist for any pair of vertices $x,y \in \Graph$, by Riesz's Lemma, as in Definition~\ref{def:vx}. 

\begin{defn}\label{def:Fin}
  For $v \in \HE$, one says that $v$ has \emph{finite support} iff there is a finite set $F \ci \Graph$ for which $v(x) = k \in \bC$ for all $x \notin F$. The set of functions of finite support in \HE is denoted $\spn\{\gd_x\}$, where $\gd_x$ is the Dirac mass at $x$, i.e., the element of \HE containing the characteristic function of the singleton $\{x\}$. It is immediate from \eqref{eqn:def:energy-form} that $\energy(\gd_x) = \cond(x)$, whence $\gd_x \in \HE$.
  Define \Fin to be the closure of $\spn\{\gd_x\}$ with respect to \energy. 
\end{defn}

\begin{defn}\label{def:Harm}
  The set of harmonic functions of finite energy is denoted
  \linenopax
  \begin{equation}\label{eqn:Harm}
    \Harm := \{v \in \HE \suth \Lap v(x) = 0, \text{ for all } x \in \Graph\}.
  \end{equation}
  It may be the case that the only harmonic functions of finite energy are constant (and hence trivial in \HE). This is true, for example, on any finite network. 
\end{defn}

\begin{lemma}[{\cite[2.11]{DGG}}]
  \label{thm:<delta_x,v>=Lapv(x)}
  For any $x \in \Graph$, one has $\la \gd_x, u \ra_\energy = \Lap u(x)$.
\end{lemma}

The following result follows easily from Lemma~\ref{thm:<delta_x,v>=Lapv(x)}; cf.~\cite[Thm.~2.15]{DGG}.

\begin{theorem}[Royden decomposition]
  \label{thm:HE=Fin+Harm}
  $\HE = \Fin \oplus \Harm$.
\end{theorem}

\begin{remark}\label{rem:dual-systems}
  The Royden decomposition illustrates one of the advantages of working with $\la u,v\ra_\energy$, as opposed to the inner product on $\ell^2(G)$ or the grounded energy product $\la u,v\ra_o := \la u,v\ra_\energy + u(o)v(o)$. Another advantage is the following: by combining \eqref{eqn:v_x} and the conclusion of Lemma~\ref{thm:<delta_x,v>=Lapv(x)}, one can reconstruct the network $(\Graph,\cond)$ (or equivalently, the corresponding Laplacian) from the dual systems (i) $(\gd_x)_{x \in X}$ and (ii) $(v_x)_{x \in X}$. Indeed, from (ii), we obtain the (relative) reproducing kernel Hilbert space \HE and from (ii), we get an associated operator $(\Lap u)(x) = \la \gd_x,u\ra_\energy$ for $u\in \HE$. 
\end{remark}

\begin{defn}\label{def:R(x)}
  Denote the (free) effective resistance from $x$ to $y$ by 
  \linenopax
  \begin{align}\label{eqn:R(x,y)}
    R(x,y) := R^F(x,y) = \energy(v_x - v_y) = \|v_x - v_y \|_\energy^2.
  \end{align}
  This quantity represents the voltage drop measured when one unit of current is passed into the network at $x$ and removed at $y$, and the central equality in \eqref{eqn:R(x,y)} is proved in \cite{ERM} and elsewhere in the literature; see \cite{Lyons,Kig03} for different formulations.
\end{defn}

The following results will be useful in the sequel; for further details, please see \cite{DGG,ERM,bdG,RBIN} and \cite{OTERN}.

\begin{lemma}[{\cite[Lem~2.23]{DGG}}]
  \label{thm:energy-kernel-is-real}
  Every $v_x$ is \bR-valued, with $v_x(y) - v_x(o) >0$ for all $y \neq o$. 
\end{lemma}

\begin{lemma}[{\cite[Lem~6.9]{bdG}}]
  \label{thm:monopoles-and-dipoles-are-bounded}
  Every $v_x$ is bounded. In particular, if we define
  \linenopax
  \begin{align}\label{eqn:sup-norm-in-HE}
    \|u\|_\iy := \sup_{x,y \in \Graph} |u(x)-u(y)|. 
  \end{align}
  then we always have $\|v_x\|_\iy \leq R(x,o)$.
\end{lemma}

\begin{lemma}[{\cite[Lem~6.8]{bdG}}]
  \label{thm:Pfin-preserves-boundedness}
  If $v \in \HE$ is bounded, then $\Pfin v$ is also bounded. 
\end{lemma}

\begin{defn}\label{def:p(x,y)}
  Let $p(x,y) := \frac{\cond_{xy}}{\cond(x)}$ so that $p(x,y)$ defines a random walk on the network, with transition probabilities weighted by the conductances. Then let
  \linenopax
  \begin{align}\label{eqn:Prob[x->y]}
    \prob[x \to y] := \prob_x(\gt_y < \gt_x^+) 
  \end{align}
  be the probability that the random walk started at $x$ reaches $y$ before returning to $x$. In \eqref{eqn:Prob[x->y]}, $\gt_z$ is the hitting time of the vertex $z$ and $\gt_z^+ := \min\{\gt_z,1\}$.
\end{defn}

\begin{cor}[{\cite[Cor.~3.13 and Cor.~3.15]{ERM}}]
  \label{thm:c(x)R(x)=Prob[x->o]}
  For any $x \neq o$, one has
  \linenopax
  \begin{equation}\label{eqn:c(x)R(x)=Prob[x->o]}
    \prob[x \to o] = \frac1{\cond(x) R(x,o)}.
  \end{equation}
\end{cor}

%% file: conductances.tex

\section{Comparing different conductance functions}
\label{sec:comparing}

Given a network $(G,c)$, we will be interested in comparing its energy space $\HE = \HEc$ and Laplace operator $\Lap = \Lc$ with those corresponding to a different conductance function $b$. 
To be clarify dependence on the conductance functions, we use scripts to distinguish between objects corresponding to different underlying conductance functions. For example, $\Lc = \Lap$ in \eqref{eqn:Lap} and $\Ec = \energy$ in \eqref{eqn:def:energy-form}, as opposed to 
  \linenopax
  \begin{equation}\label{eqn:Lap-b}
    (\Lb v)(x) :
    = \sum_{y \nbr x} b_{xy}(v(x)-v(y)).
  \end{equation}
  and
  \linenopax
  \begin{align}\label{eqn:energy-form-b}
    \energy_b(u,v)
    = \la u,v \ra_\Eb
    = \frac12 \sum_{x,y \in \Graph} b_{xy}(\cj{u}(x)-\cj{u}(y))(v(x)-v(y)),
  \end{align}
  with domain $\dom \Eb = \{u:\Graph \to \bC \suth \Eb(u)<\iy\}$. It is clear that \HEb also depends on $b$, and so too does the energy kernel $\{v_x^{(b)}\}_{x \in G}$. We will take the domains to be
  \linenopax
  \begin{align}\label{eqn:Lap-domains}
    \dom \Lb = \spn\{\vxb\}_{x \in G} 
    \q\text{ and }\q 
    \dom \Lc = \spn\{\vxc\}_{x \in G}.
  \end{align}
  
\begin{remark}\label{rem:b-connected}
  Given a network $(G,c)$ and a new conductance function $b \leq c$, it may be that $b_{xy}=0$ even though $c_{xy} > 0$, and consequently the edge structure of $(G,b)$ may be very different from $(G,c)$. \emph{However}, we will \textbf{always} make the assumption that $(G,b)$ is connected, so that Lemma~\ref{thm:connected-gives-equal-kernels} may be applied.    
\end{remark}

\begin{defn}\label{def:b}
  Let $b:\verts \times \verts \to [0,\iy)$ be a symmetric function satisfying
  \linenopax
  \begin{align*}
    b_{xy} \leq c_{xy}, \q\text{ for all } x,y \in \verts.
  \end{align*}
  In this case, we write $b \leq c$. Note that we will always assume $(G,b)$ is connected; see Remark~\ref{rem:b-connected}.
\end{defn}

\begin{lemma}\label{thm:I-contractive}
  Inclusion gives natural contractive embedding $\sI:\HEc \hookrightarrow \HEb$.
  \begin{proof}
    Since $b \leq c$, one has 
    \linenopax
    \begin{align}\label{eqn:Eb<Ec}
      \Eb(u) 
      = \frac12 \sum_{x,y \in \Graph} b_{xy}|u(x)-u(y)|^2
      \leq \frac12 \sum_{x,y \in \Graph} c_{xy}|u(x)-u(y)|^2
      = \Ec(u)  
    \end{align}
    for any function $u:G \to \bR$, and hence $\|\sI u\|_\Eb \leq \|u\|_\Ec$.
  \end{proof}
\end{lemma}

\begin{lemma}\label{thm:I(Fin)-lies-in-Fin}
  $\sI(\Finc) \hookrightarrow \Finb$ and $\sI^\ad(\Harmb) \hookrightarrow \Harmc$.
  \begin{proof}
    The first follows from the fact that $\sI(\gd_x) = \gd_x$, whence the second follows because adjoints preserve the orthocomplements (see Theorem~\ref{thm:HE=Fin+Harm}), i.e., 
    \linenopax
    \begin{align*}
      \sI^\ad\left(\Harmb\right) 
      &= \sI^\ad\left((\Finb)^\perp\right) 
      \ci \left(\Finc\right)^\perp 
      = \Harmc.
      \qedhere
    \end{align*}
  \end{proof}
\end{lemma}

Lemma~\ref{thm:connected-gives-equal-kernels} clarifies the nature of the blanket assumption that $(G,b)$ is connected; see Remark~\ref{rem:b-connected}.
\begin{lemma}\label{thm:connected-gives-equal-kernels}
  If $(G,c)$ is a network and $b \leq c$, then the following are equivalent:
  \begin{enumerate}[(i)]
    \item $(G,b)$ is connected.
    \item $\ker \Eb = \ker \Ec = \bC \one$.
    \item $\ker \sI = 0$.
  \end{enumerate}
  \begin{proof}
    To see $(i)\iff(ii)$, observe that $\Eb(u)$ is given by a sum of nonnegative terms and hence vanishes if and only if each summand does. Thus $\Eb(u)=0$ iff $u$ is locally constant. For $(ii)\implies(iii)$, note that $\sI(u)=0$ implies $\|u\|_\Eb=0$ and hence that $u$ is a constant function, whence $u=0$ in \HEb. For $(iii)\implies(ii)$, suppose $(G,b)$ is not connected, and define $u=1$ on one component and $u=0$ on the complement. Then $\|\sI(u)\|_\Eb=0$ but $u \neq 0$ in \HEc.
  \end{proof}
\end{lemma}

\begin{lemma}\label{thm:I*vxb=vxc}
  $\sI^\ad \vxb = \vxc$, and for general $u \in \HEb$, one can compute $\sI^\ad$ via 
  \linenopax
  \begin{align}\label{eqn:I*-formula}
    (\sI^\ad u)(x)-(\sI^\ad u)(y) = \frac{b_{xy}}{c_{xy}}(u(x)-u(y)).
  \end{align}
  \begin{proof}
    For $u \in \HEc \ci \HEb$, 
    \linenopax
    \begin{align*}
      \la \sI^\ad \vxb, u\ra_\Ec 
      = \la \vxb, \sI u\ra_\Eb 
      = u(x) - u(o)
      = \la \vxc, u\ra_\Ec.
    \end{align*}
    Now for $u \in \HEb$ and $v \in \HEc$, the latter claim follows from the fact that
    \linenopax
    \begin{align*}
      \la u, \sI v\ra_\Eb
      = \frac12 \sum_{x,y \in \Graph} b_{xy}(u(x)-u(y))(v(x)-v(y))
    \end{align*}
    is equal to
    \linenopax
    \begin{align*}
      \la \sI^\ad u, v\ra_\Ec 
      &= \frac12 \sum_{x,y \in \Graph} c_{xy}((\sI^\ad u)(x)-(\sI^\ad u)(y))(v(x)-v(y)).
      \qedhere
    \end{align*}
  \end{proof}
\end{lemma}

\begin{cor}\label{thm:I-injective}
  \sI is injective. 
  \begin{proof}
    $\ker \sI = \{0\}$ because $\spn\{\vxc\} = \ran \sI^\ad$ is dense in \HEc.
  \end{proof}
\end{cor}

\begin{remark}\label{rem:I-injective}
  Corollary~\ref{thm:I-injective} may appear trivial, but it is not. Suppose $H_1$ and $H_2$ are two Hilbert spaces with the same underlying vector space $V$, but different inner products for which $\|v\|_2 \leq \|v\|_1$, for all $v \in V$. Then the identity map $\gi:V \to V$ induces an embedding $H_1 \hookrightarrow H_2$ which can fail to be injective. For example, take $H_2$ to be the Hardy space $H_+(\bD)$ on the unit disk and take $H_1$ to be $u(z)H_+(\bD)$, the image of $H_2$ under the operation of multiplication by the function $u \in H^\iy(\bD)$. That is,
  \linenopax
  \begin{align*}
    H_1 = \{uh \suth h \in H_2\},
    \qq \| uh \|_1 := \|h\|_2.
  \end{align*}
  There are functions $u \in H^\iy(\bD)$ for which $\|uh\|_1 \neq 0$ and $\|uh\|_2=0$, even when $h$ is a nonzero element of $H_2$; see \cite{Sarason94} for details. 
\end{remark}

\begin{lemma}\label{thm:lap-is-almost-Green}
  If $\gd_{xy}$ is the Kronecker delta, then 
  \linenopax
  \begin{equation}\label{eqn:<vx,Lap(vy)>=gd(xy)+1}
    \la \vxb, \Lb \vyb\ra_\Eb = \gd_{xy} + 1 = \la \vxc, \Lc \vyc\ra_\Ec,
    \qq \forall x,y \in G \less\{o\}.
  \end{equation}
  \begin{proof}
    Note that  
    \linenopax
    \begin{align*}
      \la \vxb, \Lb \vyb\ra_\Eb
      &= (\Lb \vyb)(x) - (\Lb \vyb)(o) 
      = \la \gd_x, \vyb\ra_\Eb - \la \gd_o, \vyb\ra_\Eb,
    \end{align*}
    because $\gd_x \in \HEb$ and $\la \gd_x, u\ra_\Eb = \Lb u(x)$.
    Now the result follows via
    \linenopax
    \begin{align*}
      \la\gd_x, \vyb\ra_\Eb - \la\gd_o,\vyb\ra_\Eb 
      = (\gd_x(y) - \gd_x(o)) - (\gd_o(y) - \gd_o(o)) 
      = \gd_{xy} + 1,
    \end{align*}
    since $x,y \neq o$.
  \end{proof}
\end{lemma}

\begin{lemma}\label{thm:Lb=ILcI*}
  For $1 < b \leq c$, one has $\Lb = \sI \Lc \sI^\ad$.
  \begin{proof}
    Applying Lemma~\ref{thm:lap-is-almost-Green} and Lemma~\ref{thm:I*vxb=vxc},
    \linenopax
    \begin{align*}
      \la \vxb, \Lb \vyb\ra_\Eb
       &= \la \vxc, \Lc \vyc\ra_\Ec \\
       &= \la \sI^\ad \vxb, \Lc \sI^\ad \vyb\ra_\Ec \\
       &= \la \vxb, \sI \Lc \sI^\ad \vyb\ra_\Ec.
       \qedhere
    \end{align*}
  \end{proof}
\end{lemma}

Thus we have a commuting square
  \linenopax
    \begin{align}\label{eqn:commuting-square}
    \xymatrix{
    \HEc\ar[d]_{\Lc} & \HEb \ar[l]_{\sI^\ad} \ar[d]^{\Lb = \sI \Lc \sI^\ad} \\
    \HEc \ar[r]_{\sI} & \HEb
    }
  \end{align}
Note that one can recover the dipole property of \vxb from Lemma~\ref{thm:I*vxb=vxc} and Lemma~\ref{thm:Lb=ILcI*}: $\Lb \vxb = \sI\Lc\sI^\ad \vxb = \sI\Lc\vxc = \sI(\gd_x-\gd_o)=\gd_x-\gd_o$.

\begin{cor}\label{thm:spectral-invariant}
  $\sI^\ad \in \Hom(\Harmb,\Harmc)$ is a spectral invariant.
  \begin{proof}
    This is basically a restatement of Lemma~\ref{thm:I(Fin)-lies-in-Fin}. 
  \end{proof}
\end{cor}

This spectral invariant is also apparent from the formula $\Lb = \sI \Lc \sI^\ad$ of Lemma~\ref{thm:Lb=ILcI*}. While \sI is not a norm-preserving map, it is standard from spectral theory that one can write \sI in terms of its polar decomposition as $\sI = UP$ and then $\Lb = \sI \Lc \sI^\ad$ implies that a unitary equivalence is given by $\Lb = U \Lc U^\ad$.

In the case when $\dim \Harmb = \dim \Harmc = 1$, the spectral invariant of Corollary~\ref{thm:spectral-invariant} is just a number. This is computed explicitly for the geometric integers in Example~\ref{exm:geometric-integers}.

\begin{remark}[Open Question]\label{rem:open-question}
  For a fixed conductance function $b:\verts \times \verts \to [0,\iy)$, what are the closed subspaces $\sK \ci \HEb$ such that $\sK \cong \HEc$ for some conductance functions $c$ with $b \leq c$?
\end{remark}

\begin{cor}\label{thm:Lc-bounded-implies-Lb-bounded}
  If $b \leq c$ and \Lc is bounded on \HEc, then \Lb is bounded on \HEb. 
  \begin{proof}
    Lemma~\ref{thm:Lb=ILcI*} immediately implies 
    $\|\Lb\|_{\HEb \to \HEb} \leq \|\Lc\|_{\HEc \to \HEc}.$
  \end{proof}
\end{cor}

\begin{cor}\label{thm:L1-bounded-implies-Lb-bounded}
  If $c \equiv 1$ and \Lc is bounded on \HEc, then \Lb is bounded on \HEb for any bounded conductance function $b$.
  \begin{proof}
    Writing $\|b\|_\iy$ for the supremum of $b$, we have
    \linenopax
    \begin{align*}
      b_{xy} \leq \|b\|_\iy c_{xy} = \|b\|_\iy,
    \end{align*}
    so Corollary~\ref{thm:Lc-bounded-implies-Lb-bounded} applies to the network with conductances all equal to $\|b\|_\iy$.
  \end{proof}
\end{cor}

\begin{theorem}\label{thm:bounded-containment}
  Let $c$ be an arbitrary conductance function, and let \one be the conductance function which assigns a conductance of 1 to every edge. Then \HEc is contained in $\sH_{\sE^{(\one)}}$ if and only if there is an $\ge > 0$ such that $c_{xy} \geq \ge$ for all $x,y \in G$ with $c_{xy} > 0$. 
  \begin{proof}
    For the forward direction, suppose $K<\iy$ satisfies $\|u\|_{\sE^{(\one)}}^2 \leq K\|u\|_\Ec^2$, for all $u \in \HEc$. Note that $\Ec(\gd_x) = c(x)$ follows directly from \eqref{eqn:def:energy-form}, so
    \linenopax
    \begin{align*}
      c(x)
      = \|\gd_x\|^2_\Ec
      \geq \frac1K \|\gd_x\|_{\sE^{(\one)}}^2
      \geq \frac1K
    \end{align*}
    since $\|\gd_x\|_{\sE^{(\one)}} \geq 1$ by the connectedness of the network.
     
    For the converse, 
    \linenopax
    \begin{align*}
      \|u\|_{\sE^{(\one)}}^2 
      = \frac12 \sum_{x,y \in G} (u(x)-u(y))^2 
      \leq \frac12 \sum_{x,y \in G} \frac{c_{xy}}{\ge} (u(x)-u(y))^2 
      = \frac1\ge \|u\|_\Ec^2,
    \end{align*}
    so $\sI:\HEc \to \sH_{\sE^{(\one)}}$ is a bounded operator with $\|\sI\|_{\HEc \to \sH_{\sE^{(\one)}}} \leq \frac1{\sqrt\ge}$.
  \end{proof}
\end{theorem}

\begin{exm}[Horizontally connected binary tree]
  \label{exm:Lc-not-bounded}
  This example shows that the boundedness of the conductance function is not sufficient to imply boundedness of the Laplacian, and illustrates the interplay between spectral reciprocity and effective resistance (see also \cite{SRAMO}). To begin, let $(G,b)$ be the binary tree where every edge has conductance $c_{xy}=1$. Now let $(G,c)$ be the network obtained by connecting all vertices at level $k$ with an edge of conductance $c_k$ as in Figure~\ref{fig:conntree}. The resulting network is no longer a tree, but we call it the \emph{horizontally connected binary tree} for lack of a better name. Note that $b \leq c$.
  
  Suppose that $c_k=1$ for each $k$, so $c_{xy}$ is globally constant on \edges. However, $c(x) = 2^k+2$ for $x$ in level $k$, so $c(x)$ is clearly unbounded on \verts. (As usual, level $k$ consists of all vertices in $(G,b)$ for which the shortest path to $o$ contains exactly $k$ edges.) 
  Let $K_n$ be the complete graph on $n$ vertices. Using Schur complements (for example, as in \cite{ERM, OTERN} or \cite{Kig01, Kig03}), one can compute $R_{K_n}(x,y) = 2^{1-n}$ for any $x,y \in K_n$. Consequently, it is easy to see that $R_{(G,c)}^F(x,y)$ can be made arbitrarily small by choosing $x,y$ in level $k$, for sufficiently large $k$. By spectral reciprocity (see \cite{SRAMO}), this implies that \Lc is unbounded on \HEc. Thus, this network provides an example of how boundedness of $c_{xy}$ does not imply boundedness of \Lc. For an example of how boundedness of $c_{xy}$ does not imply boundedness of \gD on other spaces, see \cite{Woj07}.
  \begin{figure}
    \centering
    \includegraphics{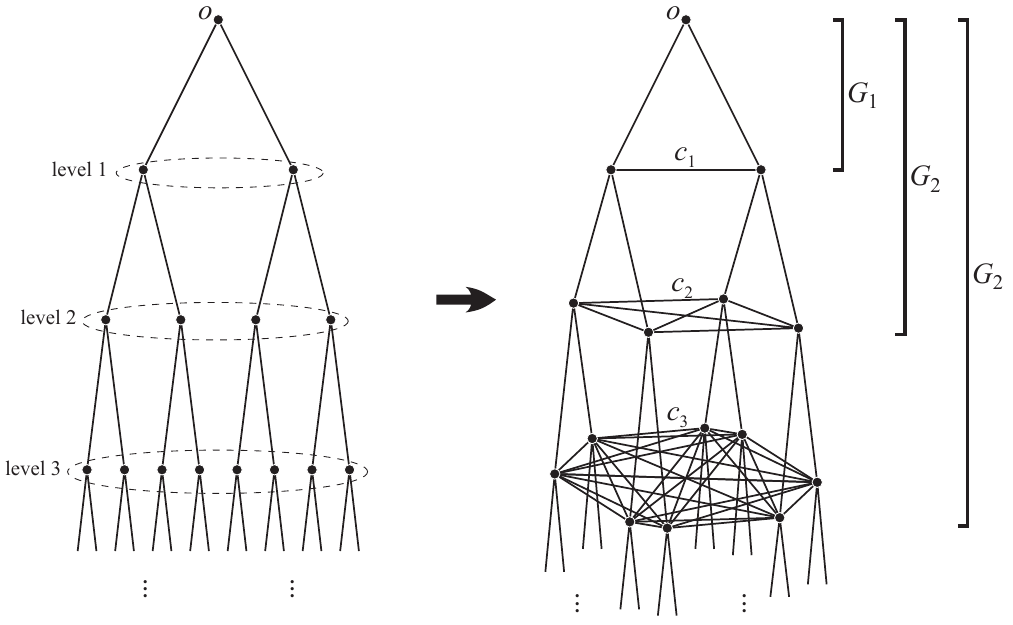}
    \caption{\captionsize Construction of the ``horizontally connected binary tree'' of Example~\ref{exm:Lc-not-bounded}. }
    \label{fig:conntree}
  \end{figure}
  
  Suppose that we choose $c_k$ so as to make $c(x)$ bounded on \verts. Then we must have $c_k = O(2^{-k})$ as $k \to \iy$, so define $c_k = 2^{-k}$. Using this, one can compute that $R_{G_k}(x,y) = 1$ for $x,y$ in level $k$ of $G_k$, for every $k$.
\end{exm}

\begin{lemma}\label{thm:self-adjointness}
  Suppose $b \leq c$. If \Lc is self-adjoint, then \Lb is self-adjoint also.
  \begin{proof}
    Take adjoints on both sides of $\Lb = \sI \Lc \sI^\ad$ (see Lemma~\ref{thm:Lb=ILcI*}). Note that the domains are as in \eqref{eqn:Lap-domains}.
  \end{proof}
\end{lemma}

\begin{exm}[Geometric integers]\label{exm:geometric-integers}
  For a fixed constant $c>1$, let $(\bZ,c^n)$ denote the network with integers for vertices, and with geometrically increasing conductances defined by $\cond_{n-1,n} = c^{\max\{|n|,|n-1|\}}$ so that the network under consideration is
  \linenopax
  \begin{align*}
    \xymatrix{
      \dots \ar@{-}[r]^{c^3}
      & -2 \ar@{-}[r]^{c^2} 
      & -1 \ar@{-}[r]^{c} 
      & 0 \ar@{-}[r]^{c} 
      & 1 \ar@{-}[r]^{c^2} 
      & 2 \ar@{-}[r]^{c^3} 
      & 3 \ar@{-}[r]^{c^4} 
      & \dots
    }
  \end{align*} 
  as in \cite[Ex.~6.2]{DGG}, and fix $o=0$. It is shown in \cite[\S4.2]{SRAMO} that \Lc is not self-adjoint, and a defect vector $\gf \in \HEc$ is constructed which satisfies
  \linenopax
  \begin{align}\label{eqn:defect}
    \Lc \gf = -\gf. 
  \end{align}
  However, for $b \equiv 1$, \Lb is bounded and Hermitian, and thus clearly self-adjoint. This example shows that the converse of Lemma~\ref{thm:self-adjointness} does not hold. Using Fourier theory, one can show that $\HEb \cong L^2\left((-\gp,\gp),\sin^2(\frac t2)\right)$; see \cite[\S6.3]{Friedrichs}, for example. 
  
  So Lemma~\ref{thm:Lb=ILcI*} gives $\Lb = \sI \Lc \sI^\ad$, where \Lb is bounded and \Lc is unbounded and not self-adjoint. The inclusion $\sI:\HEc \to \HEb$ indicates that
  \linenopax
  \begin{align*}
    \HEb = \HEc \oplus \HEc^\perp,
  \end{align*}
  where $\HEc^\perp = \HEb \ominus \HEc$, and that \Lc is a matrix corner of \Lb:
  \linenopax
  \begin{align}
    \Lb = 
    \left[\begin{array}{cc}
      \Lc & A \\ A^\ad & B 
    \end{array}\right].
  \end{align}
  This gives another way to relate the operators \Lb and \Lc.
\end{exm}

\subsection{The adjoint of \Lb with respect to \Ec}

For the results in this section we consider the adjoint of \Lb with respect to \Ec and denote it by \Lbc, in other words, we are interested in 
\linenopax
\begin{align*}
  \la \Lbc u, v \ra_{\Ec} = \la u, \Lb v \ra_{\Ec}.
\end{align*}
It will be helpful to know the action of $\sI^\ad$ on \Fin, as given in Lemma~\ref{thm:I(dirac)}; this result also generalizes the dipole property $\Lap v = \gd_x - \gd_y$ of Definition~\ref{def:dipole}.

\begin{lemma}\label{thm:I(dirac)}
  For $1 < b \leq c$, one has $\spn\{\vxc\} \ci \dom \Lbc$ and
  \linenopax
  \begin{align}\label{eqn:I(dirac)}
    \Lbc \vxc = \sI^\ad(\gd_x-\gd_o).
  \end{align}
  \begin{proof}
    For any fixed $x \in G$ and $u \in \HEc$, we have the estimate
    \linenopax
    \begin{align*}
      \la \vxc, \Lb u \ra_\Ec
      &= \Lb u(x) - \Lb u(o)
      = \la \gd_x-\gd_o, u \ra_\Eb
      \leq \|\gd_x - \gd_o\|_\Eb \cdot \|u\|_\Eb,
    \end{align*}
    by by Lemma~\ref{thm:<delta_x,v>=Lapv(x)} followed by \eqref{eqn:v_x}.
    This shows $\spn\{\vxc\} \ci \dom \Lbc$ and $\la \vxc, \Lb u \ra_\Ec = \la \gd_x-\gd_o, u \ra_\Eb$, which gives \eqref{eqn:I(dirac)}.
  \end{proof}
\end{lemma}

For Theorem~\ref{thm:Lbc-identity}, we need to define $\Lc^{-1}$ via the spectral theorem. To this end, we introduce the following blanket assumption (which remains in place for the remainder of this paper). 

\begin{assm}\label{ass:Fried}
Suppose a conductance function $c$ has been fixed. If the corresponding Laplace operator \Lc is not self-adjoint, then we replace it by the Friedrichs extension as described in \cite{Friedrichs}.
\end{assm}

With Assumption~\ref{ass:Fried} in place, we can work with \Lc as a self-adjoint operator. Then by the Spectral Theorem: for any $u \in \HEc$, there is a Borel measure $\gmc_u$ on $[0,\iy)$ such that
  \linenopax
  \begin{align}\label{eqn:spec-repn}
    \la u, \gy(\Lc)u\ra_\Ec 
    = \int_0^\iy \gy(\gl) \,d\gmc_u(\gl)
    = \int_0^\iy \gy(u) \|P(d\gl)u\|_\Ec^2,
  \end{align}
where $P$ is the projection-valued measure in the spectral resolution of \Lc. This will be useful for Theorem~\ref{thm:spectral-permanence}.
Furthermore, we also have that  
  \linenopax
  \begin{align}\label{eqn:Lc-inverse}
    \Lc^{-1} := \int_0^\iy e^{-\gl \Lc} \,d\gl.
  \end{align}
This definition of the inverse is a standard application of the spectral theorem, and is based on the fact that $\int_0^\iy e^{-\gl t} \,d\gl = \frac1t$.

\begin{theorem}\label{thm:Lbc-identity}
  For $1 < b \leq c$, one has $\Lbc = \Lc^{-1} \Lb \Lc$, where $\Lc^{-1}$ is the inverse of the Friedrichs extension, defined as in \eqref{eqn:Lc-inverse}.
  \begin{proof}
    We will first show $\Lc \Lbc = \Lb \Lc$, which is equivalent to $\sI(\Lc \Lbc - \Lb \Lc) = 0$ by Corollary~\ref{thm:I-injective}. Applying Lemma~\ref{thm:I(dirac)} and Lemma~\ref{thm:Lb=ILcI*}, one has
    \linenopax
    \begin{align*}
      \Lc \Lbc \vxc
      = \sI \Lc \Lbc \vxc 
      = \sI \Lc \sI^\ad(\gd_x-\gd_o)
      = \Lb(\gd_x-\gd_o).
    \end{align*}
    Then using the dipole property $\Lc\vxc = \gd_x-\gd_o$ yields
    \linenopax
    \begin{align*}
      \Lb(\gd_x-\gd_o)
      = \Lb(\Lc\vxc)
      = \Lb(\Lc\vxc)
      = \Lb\Lc(\vxc). 
    \end{align*}
    Now we have $\Lc \Lbc(\vxc) = \Lb \Lc(\vxc)$ for any $x$, whence $\Lc \Lbc = \Lb \Lc$ follows by the density of $\spn\{\vxc\}$ in \HEc. It follows from the preceding argument that $\Lb\Lc(\spn\{\vxc\}) \ci \dom \Lc^{-1}$, and so the proof is complete.
  \end{proof}
\end{theorem}

%% file: permanence.tex

\section{Moments of \Lc and monotonicity of spectral measures}
\label{sec:permanence}

Note that we continue to assume \Lc is a self-adjoint operator as discussed in Assumption~\ref{ass:Fried}.

\begin{lemma}\label{thm:}
  For $u=\vxc-\vyc$ and $\gy(\gl) = \gl^k$, $k=0,1,2$, we have
\linenopax
\begin{align*}
  k&=0: & \la u, u \ra_\Ec &= R^F(x,y), \\
  k&=1: & \la u, \Lc u \ra_\Ec &= 2 - 2\gd_{xy}, \\
  k&=2: & \la \vxc, \Lc^2 \vxc \ra_\Ec &= c(x) + 2 c_{xy} + c(y). 
\end{align*}
  \begin{proof}
The case $k=0$ follows immediately from \eqref{eqn:R(x,y)}. For $k=1$, \eqref{eqn:<vx,Lap(vy)>=gd(xy)+1} gives
\linenopax
\begin{align*}
  \la \vxc, &\Lc \vxc \ra_\Ec - \la \vxc, \Lc \vyc \ra_\Ec - \la \vyc, \Lc \vxc \ra_\Ec + \la \vyc, \Lc \vyc \ra_\Ec \\
  &= 2 - (\gd_{xy}+1) - (\gd_{xy}+1) + 2.
\end{align*}
For $k=2$, we use the fact that the Friedrichs extension is self-adjoint and the dipole property \eqref{eqn:v_x} to compute
\linenopax
\begin{align*}
  \la \vxc, \Lc^2 \vxc \ra_\Ec 
  &= \la \Lc \vxc, \Lc \vxc \ra_\Ec \\
  &= \la \gd_x - \gd_y, \gd_x - \gd_y \ra_\Ec  \\
  &= c(x) + 2 c_{xy} + c(y). 
\end{align*}
For the last step, we used $\energy(\gd_x) = c(x)$, which is immediate from \eqref{eqn:def:energy-form}.
  \end{proof}
\end{lemma}

In Theorem~\ref{thm:spectral-permanence}, we use $\gm_u^{(c)}$ as given by \eqref{eqn:spec-repn}. Also, let 
  \linenopax
  \begin{align}\label{eqn:kWh-moment-of-mu_u}
    m_k^{(c)}(u) := \int_0^\iy \gl^k \,d\gm_u^{(c)}
  \end{align}
be the \kth moment of $\gm_u^{(c)}$, and similarly for $\gm_u^{(b)}$. We now consider the moments of \Lap via spectral theory.

\begin{theorem}[Monotonicity of spectral measures]
  \label{thm:spectral-permanence}
  Let $(G,c)$ be a given network, and let $b \leq c$. Then
  \linenopax
  \begin{align}\label{eqn:spectral-permanence}
    m_1^{(b)}(u) = m_1^{(c)}(\sI^\ad u)
    \qq\text{and}\qq
    m_2^{(b)}(u) \leq m_2^{(c)}(\sI^\ad u).
  \end{align}
  \begin{proof}
    First, note that Lemma~\ref{thm:Lb=ILcI*} gives
    \linenopax
    \begin{align*}
      m_1^{(b)} 
      = \la u, \Lb u \ra_\Eb
      = \la u, \sI \Lc \sI^\ad u \ra_\Eb
      = \la \sI^\ad u, \Lc \sI^\ad u \ra_\Ec
      = m_1^{(c)}.
    \end{align*}
    For the second moments, using Lemma~\ref{thm:Lb=ILcI*} again gives
    \linenopax
    \begin{align*}
      m_2^{(b)} 
      = \la u, (\Lb)^2 u \ra_\Eb
      = \la u, \sI \Lc \sI^\ad \sI \Lc \sI^\ad u \ra_\Eb
      = \la \Lc^\ad \sI^\ad u, \sI^\ad \sI \Lc \sI^\ad u \ra_\Ec.
    \end{align*}
    Since $\sI^\ad\sI$ is contractive by Lemma~\ref{thm:I-contractive}, 
    \linenopax
    \begin{align*}
      \la \Lc^\ad \sI^\ad u, \sI^\ad \sI \Lc \sI^\ad u \ra_\Ec
      &\leq \|\sI^\ad\sI\| \cdot \la \Lc^\ad \sI^\ad u, \Lc \sI^\ad u \ra_\Ec \\
      &\leq \la u, \sI (\Lc)^2 \sI^\ad u \ra_\Ec,
    \end{align*}
    whence $m_2^{(b)} \leq m_2^{(c)}$.
  \end{proof}
\end{theorem}

\begin{remark}\label{rem:}
  If $b_{xy} < c_{xy}$ for some edge $(xy)$, then $m_2^{(b)}(\vxb) < m_2^{(c)}(\sI^\ad\vxb)$
\end{remark}

%% file: examples.tex

\section{Examples}
\label{sec:examples}

\begin{exm}[Geometric integers]\label{exm:geometric-integers}
  Let $(\bZ,c^n)$ be the network whose vertices are the integers with conductances given by
  \linenopax
  \begin{align*}
    c_{m,n} = \begin{cases}
      c^{\max\{|m|,|n|\}}, &|m-n|=1 \\
      0, &\text{else},
    \end{cases}
  \end{align*}
  as in the following diagram:
  \linenopax
  \begin{align*}
    \xymatrix{
      & \dots \ar@{-}[r]^{c^4}
      & \vertex{-3} \ar@{-}[r]^{c^3} 
      & \vertex{-2} \ar@{-}[r]^{c^2} 
      & \vertex{-1} \ar@{-}[r]^{c} 
      & \vertex{0} \ar@{-}[r]^{c} 
      & \vertex{1} \ar@{-}[r]^{c^2} 
      & \vertex{2} \ar@{-}[r]^{c^3} 
      & \vertex{3} \ar@{-}[r]^{c^4} 
      & \dots
    }
  \end{align*}  
  It is known that \Harm is 1-dimensional for this network; see \cite{DGG}. It was also shown in \cite{SRAMO} that \Lap is not essentially self-adjoint (as an operator on \HE) for this network. 
  
  We compare $(\bZ,b^n)$ and $(\bZ,c^n)$, where $1 < b \leq c$. In this case, $\dim \Harmb = \dim \Harmc = 1$ and we can compute the (numerical) spectral invariant of Corollary~\ref{thm:spectral-invariant}. Choose unit vectors $h_b \in \Harmb$ and $h_c \in \Harmc$: 
  \linenopax
  \begin{align}\label{eqn:norm(hb)=1}
    h_b(n) = \frac{\sgn(n)}{2\sqrt{b-1}}  \left(1-\frac1{b^{|n|}}\right),
    \qq
    h_c(n) = \frac{\sgn(n)}{2\sqrt{c-1}}  \left(1-\frac1{c^{|n|}}\right).
  \end{align}
  Now since $\la \sI^\ad h_b, u \ra_\Ec = \la h_b, u \ra_\Eb$ for all $u \in \HEc$, we have
  \linenopax
  \begin{align}\label{eqn:b-c-ansatz}
    \la h_b, \vxc[n] \ra_\Eb
    &= \la \sI^\ad h_b, \vxc[n] \ra_\Ec 
    = \la K h_c, \vxc[n] \ra_\Ec
    = K \la h_c, \vxc[n] \ra_\Ec,
  \end{align}
  following the ansatz that $\sI^\ad$ should be just a numerical constant (scaling factor). 
  Suppose for simplicity that $n>0$, as the other computation is similar. On the left side of \eqref{eqn:b-c-ansatz}, we can compute directly from \eqref{eqn:def:energy-form}:
  \linenopax
  \begin{align}
    \left\la h_b, \vxc[n] \right\ra_\Eb
    &= 2\sum_{j=1}^\iy b^j \left(\frac{1-b^{-j}}{2\sqrt{b-1}} - \frac{1-b^{1-j}}{2\sqrt{b-1}}\right) \left(\vxc[n](j)-\vxb[n](j-1)\right) \notag \\
    &= \sqrt{b-1} \vxc[n](n) 
    = \sqrt{b-1} \sum_{j=1}^n \frac1{c^n}
    = \sqrt{b-1} \frac{1-c^{-n}}{c-1}, \label{eqn:der1-left}
  \end{align}
  Meanwhile, on the right side of \eqref{eqn:b-c-ansatz}, the reproducing property gives
  \linenopax
  \begin{align}\label{eqn:der1-right}
    \la h_c, \vxc[n] \ra_\Ec
    &= h_c(n)-\vstr[2]h_c(o) 
    = \frac1{2\sqrt{c-1}} \left(1-\frac1{c^{n}}\right).
  \end{align}
  Substituting \eqref{eqn:der1-left} and \eqref{eqn:der1-right} into \eqref{eqn:b-c-ansatz} gives
  \linenopax
  \begin{align*}
    \sqrt{b-1} \frac{1-c^{-n}}{c-1}  
    &= K \frac1{2\sqrt{c-1}} \left(1-\frac1{c^{n}}\right), 
  \end{align*}
  and so the corresponding spectral invariant is 
  \linenopax
  \begin{align*}
    K = \left\| \restr{\sI^\ad}{\Harmb} \,\vstr[2.2]\right\| 
    = \sqrt{\frac{1-b}{1-c}},
  \end{align*}
  and this is the factor by which $\sI^\ad$ scales the basis vector $h_b$; see Corollary~\ref{thm:spectral-invariant}.
\end{exm}